\documentstyle[12pt,amstex,amscd,graphicx]{article}

\begin{document}

\newtheorem{theorem}{Theorem}[section]
\newtheorem{prop}[theorem]{Proposition}
\newtheorem{lemma}[theorem]{Lemma}
\newtheorem{cor}[theorem]{Corollary}

\title{Cannon-Thurston Maps and Bounded Geometry}

\author{Mahan Mj}\date{}

\maketitle

\begin{abstract} 
This is an expository paper. We prove the Cannon-Thurston property for
bounded geometry surface groups with or without punctures. We  prove
three theorems, due to Cannon-Thurston, Minsky and Bowditch. The proofs
are culled out of earlier work of the author.
\end{abstract}

\begin{center}
AMS subject classification =   20F32, 57M50
\end{center}

\tableofcontents

\section{Introduction}

Let $M$ be a closed hyperbolic 3-manifold fibering over the circle with 
fiber $F$. Let $\widetilde F$ and $\widetilde M$ denote the universal
covers of $F$ and $M$ respectively. Then $\widetilde F$ and $\widetilde M$
are quasi-isometric to ${\Bbb{H}}^2$ and ${\Bbb{H}}^3$ respectively. Now let
${{\Bbb{D}}^2}={\Bbb{H}}^2\cup{\Bbb{S}}^1_\infty$ and 
${{\Bbb{D}}^3}={\Bbb{H}}^3\cup{\Bbb{S}}^2_\infty$
denote the standard compactifications. In \cite{CT} Cannon and Thurston
show that the usual inclusion of $\widetilde F$ into $\widetilde M$
extends to a continuous map from ${\Bbb{D}}^2$ to ${\Bbb{D}}^3$.
This was extended to Kleinian surface groups of bounded geometry
without parabolics by Minsky \cite{minsky-jams}. Bowditch
\cite{bowditch-ct} \cite{bowditch-stacks} proved the Cannon-Thurston
property for bounded geometry surface groups with parabolics.

In \cite{mitra-ct}, \cite{mitra-trees}, \cite{brahma-pared},
\cite{brahma-ibdd}, \cite{brahma-amalgeo}, we have given a different
approach to the Cannon-Thurston problem. Though   the
theorems of Cannon-Thurston, Minsky and Bowditch can be deduced from
ours, it might be instructive to write down a complete proof of these
results. In some ways, the proof here is easier and more
minimalistic. Another reason for writing this paper is that Cannon and
Thurston's original result \cite{CT} is unpublished. It seems only
appropriate that the theorem  that motivated all the above results be
available. 

Much of what follows is true in the setting of hyperbolic metric
spaces in the sense of Gromov. We shall often state results in this
generality.

\section{Preliminaries}

\subsection{Hyperbolic Metric Spaces}

We start off with some preliminaries about hyperbolic metric
spaces  in the sense
of Gromov \cite{gromov-hypgps}. For details, see \cite{CDP}, \cite{GhH}. Let $(X,d)$
be a hyperbolic metric space. The 
{\bf Gromov boundary} of 
 $X$, denoted by $\partial{X}$,
is the collection of equivalence classes of geodesic rays $r:[0,\infty)
\rightarrow X$ with $r(0)=x_0$ for some fixed ${x_0}\in{X}$,
where rays $r_1$
and $r_2$ are equivalent if $sup\{ d(r_1(t),r_2(t))\}<\infty$.
Let $\widehat{X}$=$X\cup\partial{X}$ denote the natural
 compactification of $X$ topologized the usual way(cf.\cite{GhH} pg. 124).

The {\bf Gromov inner product}
 of elements $a$ and $b$ relative to $c$ is defined 
by 
\begin{center}
$(a,b)_c$=1/2$[d(a,c)+d(b,c)-d(a,b)]$.
\end{center}

\medskip

{\bf Definitions:} {\it A subset $Z$ of $X$ is said to be 
{\bf $k$-quasiconvex}
if any geodesic joining $a,b\in Z$ lies in a $k$-neighborhood of $Z$.
A subset $Z$ is {\bf quasiconvex} if it is $k$-quasiconvex for some $k$.
A map $f$ from one metric space $(Y,{d_Y})$ into another metric space 
$(Z,{d_Z})$ is said to be
 a {\bf $(K,\epsilon)$-quasi-isometric embedding} if
 
\begin{center}
${\frac{1}{K}}({d_Y}({y_1},{y_2}))-\epsilon\leq{d_Z}(f({y_1}),f({y_2}))\leq{K}{d_Y}({y_1},{y_2})+\epsilon$
\end{center}
If  $f$ is a quasi-isometric embedding, 
 and every point of $Z$ lies at a uniformly bounded distance
from some $f(y)$ then $f$ is said to be a {\bf quasi-isometry}.
A $(K,{\epsilon})$-quasi-isometric embedding that is a quasi-isometry
will be called a $(K,{\epsilon})$-quasi-isometry.

A {\bf $(K,\epsilon)$-quasigeodesic}
 is a $(K,\epsilon)$-quasi-isometric embedding
of
a closed interval in $\Bbb{R}$. A $(K,0)$-quasigeodesic will also be called
a $K$-quasigeodesic.}

\medskip

Let $(X,{d_X})$ be a hyperbolic metric space and $Y$ be a subspace that is
hyperbolic with the inherited path metric $d_Y$.
By 
adjoining the Gromov boundaries $\partial{X}$ and $\partial{Y}$
 to $X$ and $Y$, one obtains their compactifications
$\widehat{X}$ and $\widehat{Y}$ respectively.

Let $ i :Y \rightarrow X$ denote inclusion.

\medskip

{\bf Definition:} { \it Let $X$ and $Y$ be hyperbolic metric spaces and
$i : Y \rightarrow X$ be an embedding. 
 A {\bf Cannon-Thurston map} $\hat{i}$  from $\widehat{Y}$ to
 $\widehat{X}$ is a continuous extension of $i$. }

\medskip

The following  lemma
 says that a Cannon-Thurston map exists
if for all $M > 0$ and $y{\in}Y$, there exists $N > 0$ such that if $\lambda$
lies outside an $N$ ball around $y$ in $Y$ then
any geodesic in $X$ joining the end-points of $\lambda$ lies
outside the $M$ ball around $i(y)$ in $X$.
For convenience of use later on, we state this somewhat
differently. The proof is similar to Lemma 2.1 of \cite{mitra-ct}.

\begin{lemma}
A Cannon-Thurston map from $\widehat{Y}$ to $\widehat{X}$
 exists if  the following condition is satisfied:

Given ${y_0}\in{Y}$, there exists a non-negative function  $M(N)$, such that 
 $M(N)\rightarrow\infty$ as $N\rightarrow\infty$ and for all geodesic segments
 $\lambda$  lying outside an $N$-ball
around ${y_0}\in{Y}$  any geodesic segment in $X$ joining
the end-points of $i(\lambda)$ lies outside the $M(N)$-ball around 
$i({y_0})\in{X}$.

\label{contlemma}
\end{lemma}

\noindent {\it Proof:}
 Suppose $i:Y\rightarrow{X}$
does not extend continuously . Since $i$ is proper, there exist 
sequences $x_m$, $y_m$ $\in{Y}$ and $p\in\partial{Y}$,
such that $x_m\rightarrow p$
and $y_m\rightarrow p$ in $\widehat{Y}$, but $i(x_m)\rightarrow u$
and $i(y_m)\rightarrow v$ in $\widehat{X}$, where 
$u,v\in\partial{X}$ and $u\neq v$.

Since $x_m\rightarrow p$ and $y_m\rightarrow p$, any geodesic in $Y$
joining $x_m$ and $y_m$ lies outside an $N_m$-ball ${y_0}\in{Y}$,
 where $N_m\rightarrow\infty$ as $m\rightarrow\infty$. Any
bi-infinite geodesic in $X$ joining  $u,v\in\partial{X}$
has to pass through some $M$-ball around $i({y_0})$ in $X$ as
$u\neq v$. There exist constants $c$ and $L$ such that for all $m > L$
any geodesic joining $i(x_m)$ and $i(y_m)$ in $X$ 
passes through an $(M+c)$-neighborhood
 of $i({y_0})$. 
Since $(M+c)$ is a constant not depending on the index $m$
this proves the lemma. $\Box$

\medskip

The above result can be interpreted as saying that a Cannon-Thurston map 
exists if the space of geodesic segments in $Y$ embeds properly in the
space of geodesic segments in $X$.

\smallskip

We shall be needing the notion of a tree of hyperbolic metric spaces \cite{BF}.

\smallskip

\noindent
{\bf Definition:}  A  tree $T$ of hyperbolic metric spaces satisfying
the {\bf q(uasi) i(sometrically) embedded condition} is a metric space $(X,d)$
admitting a map $P : X \rightarrow T$ onto a simplicial tree $T$, such
that there exist $\delta{,} \epsilon$ and $K > 0$ satisfying the following: \\
\begin{enumerate}
\item  For all vertices $v\in{T}$, 
$X_v = P^{-1}(v) \subset X$ with the induced path metric $d_v$ is a 
$\delta$-hyperbolic metric space. Further, the
inclusions ${i_v}:{X_v}\rightarrow{X}$ 
are uniformly proper, i.e. for all $M > 0$, 
there exists $N > 0$ such that for all
$v\in{T}$ and $x, y\in{X_v}$,
 $d({i_v}(x),{i_v}(y)) \leq M$ implies
${d_v}(x,y) \leq N$.
\item Let $e$ be an edge of $T$ with initial and final vertices $v_1$ and
$v_2$ respectively.
Let $X_e$ be the pre-image under $P$ of the mid-point of  $e$.  
Then $X_e$ with the induced path metric is $\delta$-hyperbolic.
\item There exist maps ${f_e}:{X_e}{\times}[0,1]\rightarrow{X}$, such that
$f_e{|}_{{X_e}{\times}(0,1)}$ is an isometry onto the pre-image of the
interior of $e$ equipped with the path metric.
\item ${f_e}|_{{X_e}{\times}\{{0}\}}$ and 
${f_e}|_{{X_e}{\times}\{{1}\}}$ are $(K,{\epsilon})$-quasi-isometric
embeddings into $X_{v_1}$ and $X_{v_2}$ respectively.
${f_e}|_{{X_e}{\times}\{{0}\}}$ and 
${f_e}|_{{X_e}{\times}\{{1}\}}$ will occasionally be referred to as
$f_{v_1}$ and $f_{v_2}$ respectively.
\end{enumerate}   

$d_v$ and $d_e$ will denote path metrics on $X_v$ and $X_e$ respectively.
$i_v$, $i_e$ will denote inclusion of $X_v$, $X_e$ respectively into $X$.

A few general lemmas  about hyperbolic metric
spaces will be useful. We reproduce the proofs here from \cite{mitra-trees}.

\smallskip

{\bf Nearest Point Projections}

\smallskip

The following Lemma says nearest point projections in a $\delta$-hyperbolic
metric space do not increase distances much.
It is a standard fact that any geodesic metric space
is quasi-isometric to a graph \cite{br-ha}. In what follows, we shall often implicitly identify spaces
with their graph approximations. Further, graphs are
declared to have edge length one. Also, without loss of
generality, we may look at functions restricted to the
vertex sets of these graphs.

\begin{lemma}
Let $(Y,d)$ be a $\delta$-hyperbolic metric space
 and  let $\mu\subset{Y}$ be
 a geodesic segment.
Let ${\pi}:Y\rightarrow\mu$ map $y\in{Y}$ to a point on
$\mu$ nearest to $y$. Then $d{(\pi{(x)},\pi{(y)})}\leq{C_3}d{(x,y)}$ for
all $x, y\in{Y}$ where $C_3$ depends only on $\delta$.
\label{easyprojnlemma}
\end{lemma}

{\it Proof:}
Let $[a,b]$ denote a geodesic edge-path joining vertices $a, b$. Recall that
the Gromov inner product $(a,b{)}_c$=1/2[$d{(a,c)}+{d}(b,c)-{d}(a,b)]$.
It suffices by repeated use of the triangle inequality to prove the Lemma
when ${d}(x,y) \leq 1$. 
Let $u, v, w$ be points on $[x,\pi{(x)}]$, $[{\pi}(x),{\pi}(y)]$ and 
$[{\pi}(y),x]$ respectively such that ${d}(u,{\pi}(x))={d}(v,{\pi}(x))$,
${d}(v,{\pi}(y))={d}(w,{\pi}(y))$ and ${d}(w,x)={d}(u,x)$.
Then ${(x,{\pi}(y))}_{\pi{(x)}}={d}(u,{\pi}(x))$. Also, since $Y$
 is $\delta$-hyperbolic, 
the diameter of the inscribed triangle with vertices $u, v, w$ is
less than or equal to $2\delta$ (See \cite{Shortetal}).

\begin{eqnarray*}
{d}(u,x)+{d}(u,v) & \geq & {d}(x,{\pi}(x)) =  {d}(u,x)+{d}(u,{\pi}(x)) \\
\Rightarrow {d}(u,{\pi}(x)) & \leq & {d}(u,v)\leq{2\delta}  \\
\Rightarrow {(x,{\pi}(y))}_{{\pi}(x)} & \leq & 2{\delta}
\end{eqnarray*}

Similarly, ${(y,{\pi}(x))}_{\pi{(y)}}\leq{2\delta}$.

\begin{center}
i.e. ${d}(x,{\pi}(x))+{d}({\pi}(x),{\pi}(y))-{d}(x,{\pi}(y))\leq{4\delta}$

and  ${d}(y,{\pi}(y))+{d}({\pi}(x),{\pi}(y))-{d}(y,{\pi}(x))\leq{4\delta}$
\end{center}

Therefore,
\begin{eqnarray*}
\lefteqn{2{d}({\pi}(x),{\pi}(y)) }   \\
    & \leq & {8\delta}+{d}(x,{\pi}(y))-{d}(y,{\pi}(y))+{d}(y,{\pi}(x))-{d}(x,{\pi}(x)) \\
    & \leq & {8\delta}+{d}(x,y)+{d}(x,y) \\
     & \leq & {8\delta}+2     
\end{eqnarray*}
 Hence  ${d}({\pi}(x),{\pi}(y))\leq{4\delta}+1$.
Choosing $C_3 = {4\delta}+1$, we are through. $\Box$

\medskip

\begin{lemma}
Let $(Y,d)$ be a $\delta$-hyperbolic metric space.
Let $\mu$ be a geodesic segment in $Y$ with end-points $a, b$ and let
$x$ be any vertex in $Y$. Let $y$ be a vertex on $\mu$ such that 
${d}(x,y)\leq{d}(x,z)$ for any $z\in\mu$. Then a geodesic path
from $x$ to $y$ followed by a geodesic path from $y$ to $z$ is a 
$k$-quasigeodesic for some $k$ dependent only on $\delta$.
\label{unionofgeodslemma}
\end{lemma}

\noindent {\bf Proof:}
As in Lemma \ref{easyprojnlemma}, let $u, v, w$ be points on edges $[x,y]$,
$[y,z]$ and $[z,x]$ respectively such that ${d}(u,y)={d}(v,y)$,
 ${d}(v,z)={d}(w,z)$ and  ${d}(w,x)={d}(u,x)$.
Then ${d}(u,y)={(z,x)}_y\leq{2\delta}$ and the inscribed triangle with vertices
$u, v, w$ has diameter less than or equal to $2\delta$ (See \cite{Shortetal}).
 $[x,y]\cup{[y,z]}$ is a union of 2 geodesic paths lying in a 
$4\delta$ neighborhood of a geodesic $[x,z]$. Hence a geodesic path from
$x$ to $y$ followed by a geodesic path from $y$ to $z$ is a $k-$quasigeodesic
for some $k$ dependent only on $\delta$. $\Box$

\medskip

\begin{lemma}
Suppose $(Y,d)$ is a $\delta$-hyperbolic metric space.
If $\mu$ is a $({k_0},{\epsilon_0})$-quasigeodesic in $Y$ and $p, q, r$ are 
3 points in order on $\mu$ then ${(p,r)}_q\leq{k_1}$ for some $k_1$
dependent on $k_0$, $\epsilon_0$ and $\delta$ only.
\label{qgeodiplemma}
\end{lemma}

\noindent {\bf Proof:}  $[a,b]$ will denote a geodesic path joining $a, b$. 
Since $p, q, r$ are 3 points in order on $\mu$, $[p,q]$ followed
by $[q,r]$ is a $({k^{\prime}},{\epsilon^{\prime}})$-quasigeodesic in the $\delta$-hyperbolic metric space $Y$ where $({k^{\prime}},{\epsilon^{\prime}})$ depend only
on
$ \delta , {k_0}, {\epsilon_0}$. 
Hence there exists a $k_1$ dependent on $k_0$, $\epsilon_0$
 and $\delta$ alone
such that ${d}(q,[p,r])\leq{k_1}$. Let $s$ be a point on $[p,r]$ such that
${d}(q,s)={d}(q,[p,r])\leq{k_1}$. Then 
\begin{eqnarray*}
{(p,r)}_q & = & 1/2({d}(p,q)+{d}(r,q)-{d}(p,r))  \\
           & = & 1/2({d}(p,q)+{d}(r,q)-{d}(p,s)-{d}(r,s)) \\
         & \leq & {d}(q,s)\leq{k_1}. \Box
\end{eqnarray*}

\subsection{Stability of Tripods, or NPP's and QI's almost commute}

A crucial property of hyperbolic metric spaces is {\em stability of
  quasigeodesics}, i.e. any quasigeodesic (in particular an image of a
  geodesic under a quasi-isometry) lies in a bounded neighborhood of a
  geodesic. This property can be easily extended to quasiconvex sets
  \cite{GhH} \cite{CDP}. Here we are interested in a particular kind
  of quasiconvex set, the {\bf tripod}. In general, a
 tripod is a union of three
  geodesic segments, all of which share a common end-point. It is easy
  to see (for instance by thinness of triangles) that a tripod is
  quasiconvex. We shall be interested in a special kind of a
  tripod. Let $[a,b]$ be a geodesic in a hyperbolic metric space
  $X$. Let $x \in X$ be some point. Let $p$ be a nearest point
  projection of $x$ onto $[a,b]$. We shall look at tripods of the form
  $[a,b]\cup [x,p]$. We shall show that such tripods are stable under
  quasi-isometries. 

However, we shall interpret this differently to say that nearest point
projections and quasi-isometries 
in hyperbolic metric spaces `almost commute'.  
The following Lemma says precisely this:
 nearest point projections and quasi-isometries
in hyperbolic metric spaces `almost commute'. 

\begin{lemma}
Suppose $(Y,d)$ is $\delta$-hyperbolic.
Let $\mu_1$ be some geodesic segment in $Y$ joining $a, b$ and let $p$
be any vertex of $Y$. Also let $q$ be a vertex on $\mu_1$ such that
${d}(p,q)\leq{d}(p,x)$ for $x\in\mu_1$. 
Let $\phi$ be a $(K,{\epsilon})$ - quasi-isometry from $Y$ to itself.
Let $\mu_2$ be a geodesic segment 
in $Y$ joining ${\phi}(a)$ to ${\phi}(b)$ for some $g\in{S}$. Let
$r$ be a point on $\mu_2$ such that ${d}({\phi}(p),r)\leq{d}({\phi(p)},x)$ for $x\in\mu_2$.
Then ${d}(r,{\phi}(q))\leq{C_4}$ for some constant $C_4$ depending   only on
$K, \epsilon $ and $\delta$. 
\label{cruciallemma}
\end{lemma}

\noindent {\bf Proof:}
Since  ${\phi}({\mu_1})$ is a 
$({K,\epsilon})-$   quasigeodesic
joining ${\phi}(a)$ to ${\phi}(b)$, it lies in a $K^{\prime}$-neighborhood 
of $\mu_2$ where $K^{\prime}$ depends only on $K, {\epsilon}, \delta$. 
Let $u$ be a point
in ${\phi}({\mu_1})$ lying at a distance at most $K^{\prime}$ from $r$.
Without loss of generality suppose that $u$ lies on ${\phi}([q,b])$, 
where $[q,b]$ denotes the geodesic subsegment of $\mu_1$ joining $q, b$.
[See Figure 1 below.]

\smallskip

\begin{center}

\includegraphics[height=5cm]{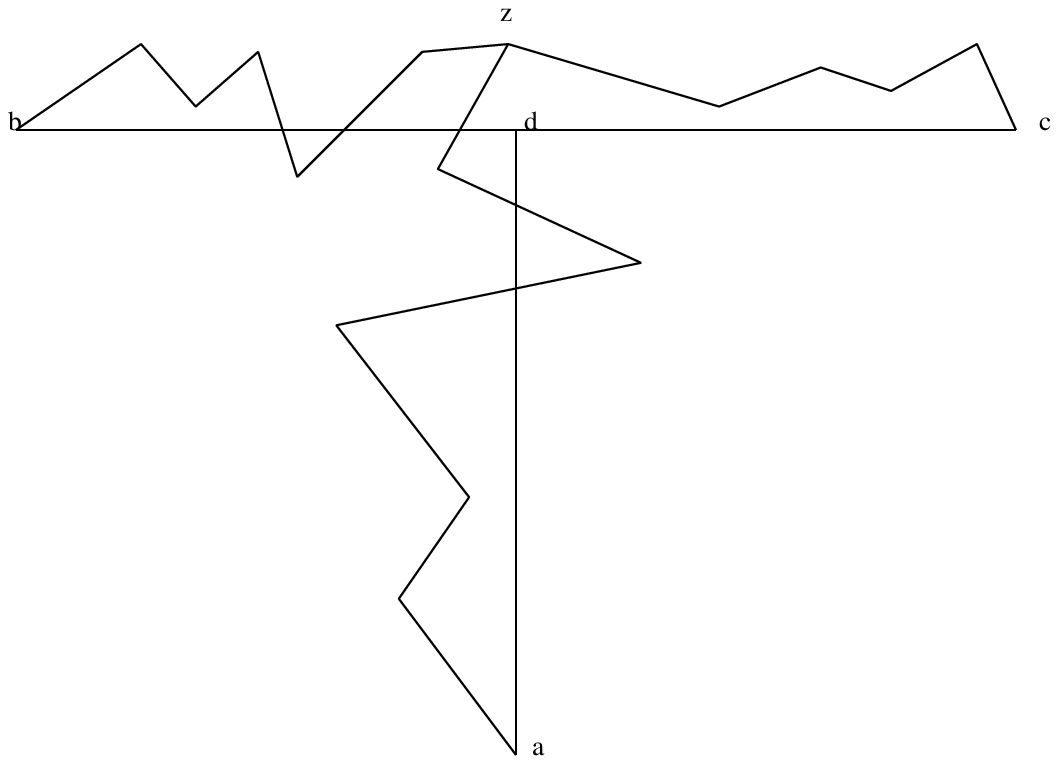}

\smallskip

{\em Figure 1: Quasi-isometries and Nearest Point Projections Almost Commute}

\end{center}

\smallskip

Let $[p,q]$ denote a geodesic joining $p, q$.
From Lemma \ref{unionofgeodslemma} $[p,q]\cup{[q,b]}$ is a $k$-quasigeodesic,
where $k$ depends on $\delta$ alone. Therefore 
${\phi}([p,q])\cup{\phi}([q,b])$ is a 
$({{K_0}, \epsilon_0})$-quasigeodesic, where $K_0, {\epsilon_0}$ depend
on $K, k, \epsilon$. Hence, by Lemma \ref{qgeodiplemma} 
${({\phi}(p),u)}_{{\phi}(q)}\leq{K_1}$, where $K_1$ depends on $K$, $k$,
$\epsilon$ and $\delta$ alone. Therefore, 
\begin{eqnarray*}
\lefteqn{ {({\phi}(p),r)}_{{\phi}(q)} } \\
  & = & 1/2[{d}({\phi}(p),{\phi}(q))+{d}(r,{\phi}(q))-{d}(r,{\phi}(p))]  \\
  & \leq & 1/2[{d}({\phi}(p),{\phi}(q))+{d}(u,{\phi}(q))+{d}(r,u) \\ 
  &      &  \hspace{1.5in}       -{d}(u,{\phi}(p))+{d}(r,u)] \\
  & = & {({\phi}(p),u)}_{{\phi}(q)}+{d}(r,u) \\
  & \leq & {K_1}+{K^{\prime}}
\end{eqnarray*}

There exists $s\in{\mu_2}$ such that ${d}(s,{\phi}(q))\leq{K^{\prime}}$

\begin{eqnarray*}
{{({\phi}(p),r)}_{s}} & = & 1/2[{d}({\phi}(p),s)+{d}(r,s)-{d}(r,{\phi}(p))] \\
   & \leq &  1/2[{d}({\phi}(p),{\phi}(q))+{d}(r,{\phi}(q))-{d}(r,{\phi}(p))]+{K^{\prime}}  \\
  & = & {({\phi}(p),r)}_{{\phi}(q)}+{K^{\prime}} \\
  & \leq &  {K_1}+{K^{\prime}}+{K^{\prime}} \\
  & = &  {K_1}+2{K^{\prime}}
\end{eqnarray*}

Also, as in the proof of Lemma \ref{easyprojnlemma} 
${({\phi}(p),s)}_r\leq{2\delta}$

\begin{eqnarray*}
{d}(r,s) & = & {({\phi}(p),s)}_{r}{+}{({\phi}(p),r)}_{s}  \\
             & \leq & K_1{+}2{K^{\prime}}{+}2{\delta} \\
{d}(r,{\phi}(q)) & \leq & K_1{+}2{K^{\prime}}{+}2{\delta}{+}{d}(s,{\phi}(q)) \\
                         & \leq & K_1{+}2{K}^{\prime}{+}2{\delta}{+}{K}^{\prime}  
\end{eqnarray*}

Let $C_4{=}K_1{+}3{K^{\prime}}{+}2{\delta}$. Then ${d}(r,{\phi}(q))\leq{C_4}$ 
and $C_4$ is independent of $a, b, p$. $\Box$

\section{Rays  of Spaces and Hyperbolic Ladders}

\subsection{Trees of Hyperbolic Metric spaces}

In this section, we shall consider trees $T$ of hyperbolic metric
spaces satisfying the qi-embedded condition. $(X,d_X)$ will denote the space
that is a tree of hyperbolic metric spaces. 
Our trees will be rather special. 

\smallskip

\noindent $\bullet$ {\bf $T$ will  either be a ray
$[0, \infty )$ or a bi-infinite geodesic $(- \infty , \infty )$.}

\smallskip

 The
  vertices of $T$ will be the integer points $j \in \{ 0 \} \cup
 {\Bbb{N}}$ or
  $\Bbb{Z}$ according as $T$ is 
$[0, \infty )$ or  $(- \infty , \infty )$. The edges will be of the
   form $[j, j+1]$. All the vertex and edge spaces will be abstractly
   isometric and identified with a fixed hyperbolic metric space
   $Y$. The vertex space over $j$ will be denoted as $Y_{j}$, and the
   edge space over $[j, j+1]$ will be denoted as $Y_{ej}$. There are
   two qi-embeddings of $Y_{ej}$, one into $Y_j$ (denoted $i_{j0}$), and
   the other into 
   $Y_{i+1}$ (denoted $i_{j1}$). We shall demand:

\begin{enumerate}
\item $i_{j0}$ is the {\bf identity} map for  all $j$, 
i.e the identification
  of $Y_j$ and $Y_{ej}$ with $Y$ is the same \\
\item There exist $K, \epsilon$ such that for all $j $, $i_{j1}$ is a
  $(K, \epsilon)$ quasi-isometry between $Y_{ej}$ and $Y_{j+1}$. \\
\end{enumerate}

Thus, the tree $T$ is assumed rooted at $0$, and maps are described on
this basis.

This induces a map $\phi_j$ from $Y_j$ to $Y_{j+1}$ which is a uniform
$(K, \epsilon)$ quasi-isometry for all $j $.
  Let $\Phi_j$ denote the
induced map on geodesics. So $\Phi_j([a,b])$ is a geodesic in
$Y_{j+1}$ joining $\phi_j (a), \phi_{j+1} (b)$. 
$\phi_j^{-1}$ and $\Phi_j^{-1}$ will denote the quasi-isometric
inverse of $\phi_j$ and the induced map on geodesics respectively. 
We shall assume that the quasi-isometric inverse $\phi_j^{-1}$ is also
a $(K, \epsilon )$ quasi-isometry. {\it (Note that 
$\phi_j^{-1}$ is just a notation for the 
quasi-isometric inverse of $\phi$ and is not necessarily 
a set-theoretic inverse.)}

\smallskip

\noindent $\bullet$ The space $Y$ in question will also be quite
special. {\bf $Y$ will either be the universal cover of a closed hyperbolic
surface (hence the hyperbolic plane), or the universal cover of a
finite volume hyperbolic surface minus  cusps. }

\smallskip

The first case, where $Y$ is ${\Bbb{H}}^2$ will be necessary when we
prove the results of Cannon-Thurston\cite{CT} and Minsky
\cite{minsky-jams}. The second case, where $Y$ is  ${\Bbb{H}}^2$ minus
an equivariant collection of horodisks, will be useful while proving
a result of Bowditch \cite{bowditch-ct}.

\subsection{Hyperbolic Ladders}

Given a geodesic $\lambda = \lambda_0 \subset Y_0$, we shall now construct a
hyperbolic ladder $B_\lambda \subset X$ containing $\lambda_0$. We
shall then prove that $B_\lambda$ is uniformly quasiconvex
(independent of $\lambda$). 
Inductively define:

\begin{center}
$\lambda_{j+1} = \Phi_j (\lambda_j )$ for $j \geq 0$ \\
$\lambda_{j-1} = \Phi_j^{-1} (\lambda_j )$ for $j \leq 0$ \\
$B_\lambda = \bigcup_j \lambda_j$
\end{center}

[See Figure 2 below.]

\smallskip

\begin{center}

\includegraphics[height=5cm]{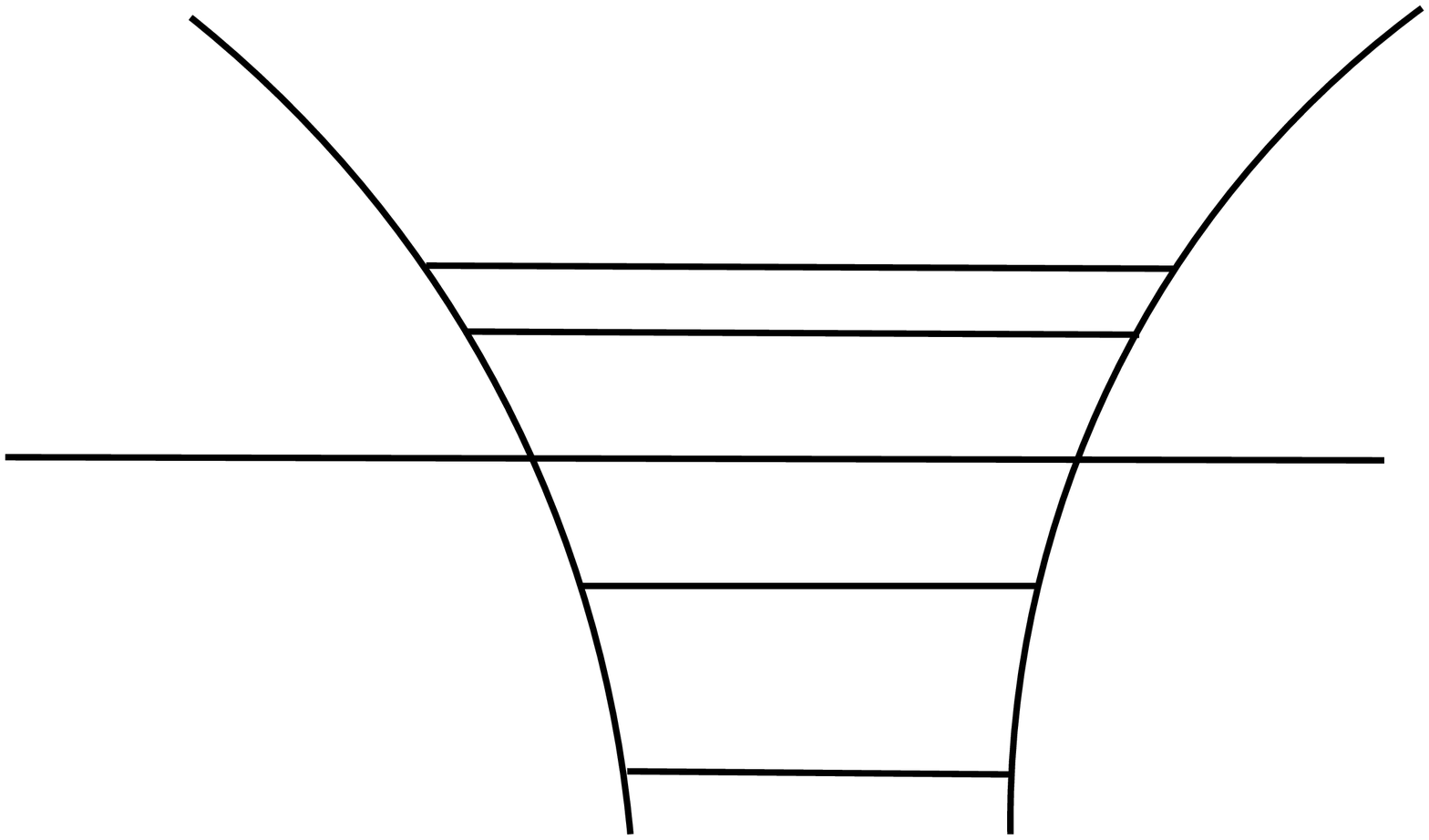}

\smallskip

{\em Figure 2: The Hyperbolic Ladder}

\end{center}

\smallskip

To prove quasiconvexity, we construct a retraction $\Pi_\lambda$ from
$\bigcup_j Y_j$ onto $B_\lambda$ that fixes $B_\lambda$ and does not
stretch distances much. So $\Pi_\lambda$ is a `quasi-Lipschitz' map.

On $Y_j$ define $\pi_j (y)$ to be a nearest-point projection of $y$
onto $\lambda_j$. Thus, $d_{Y_j}(y,\pi_j(y)) \leq 
d_{Y_j}(y,z)$ for all $z \in
\lambda_j$. 

Next, define

\begin{center}

$\Pi_\lambda (y) = \pi_j (y) $ for $y \in Y_j$.

\end{center}

We shall choose a collection of {\bf admissible paths} in $X$, such
that any path in $X$ can be uniformly approximated by these. An
admissible path in $X$ is a union of paths of the following forms:\\

\begin{enumerate}

\item Paths lying in $Y_j$ for some $j$. \\
\item Paths of length one joining $y \in Y_j$ to $\phi_j (y) \in
  Y_{j+1}$ for $j \geq 0$ \\
\item Paths of length one joining $y \in Y_j$ to $\phi_j^{-1} (y) \in
  Y_{j-1}$ for $j \leq 0$ \\
\end{enumerate}

\noindent {\bf Note:} Any path in $X $ can be uniformly approximated
by an admissible path. 

\noindent {\bf Definition:} An admissible geodesic is an admissible
path that minimizes distances amongst admissible path. 

\smallskip

We are now in a position to prove the main technical theorem of this
paper.

\begin{theorem}
There exists $C_0\geq{0}$ such that 
${d}_X({\Pi_\lambda}(x),{\Pi_\lambda}(y))\leq C_0{d_X}(x,y)$ for $x, y
\in \bigcup_i Y_i$.
Further, given $\delta \geq 0$, there exists $C \geq 0$ such that if
$X$ is $\delta$-hyperbolic, then $B_\lambda$ is $C$-quasiconvex.
\label{mainref}
\end{theorem}

\noindent {\bf Proof:} It suffices to prove the 
the theorem when ${d_X}(x,y) \leq 1$.  

\medskip

{\it Case (a):} $x, y \in Y_j$ for some $j$: \\
 From Lemma \ref{easyprojnlemma}, 
there exists $C_3$ such that
${d_{Y_j}}(\pi_j(x),\pi_j(y))\leq{C_3}$.
Since embeddings of $Y_j$ in $X$ are 1-Lipschitz (edges go to edges), 
$d_X({\Pi_{\lambda}}(x),{\Pi_{\lambda}}(y)) \leq C_3$.
	
{\it Case (b):} $x\in Y_j$ and $y = \phi_j (x) \in Y_{j+1}$
for some $j \geq 0$: \\
By Lemma \ref{cruciallemma}, there exists $C_4 \geq 0$ such that \\
\begin{center}
$d_{Y}(\phi_j (\pi_j (x)), {\pi_{j+1}}( \phi_j(x))) \leq  C_4$
\end{center}
 Unravelling definitions, and using the uniformity of the embedding
of $Y_j$ in $X$, \\
 $d({{\Pi_{\lambda}}}(x),{{\Pi_{\lambda}}}(y)) \leq C_4 + 1$. \\

{\it Case (c):} $x\in Y_j$ and $y = \phi_j^{-1} (x) \in Y_{j-1}$
for some $j \leq 0$: \\
The argument is just as in Case(b) above. 

\smallskip

Finally, to prove quasiconvexity of $B_\lambda$, we consider any two
points $a, b \in B_\lambda$. Let $\mu$ be a geodesic in $X$ joining
$a, b$. Then by the above argument, $\Pi_\lambda (\mu \cap \bigcup_j
Y_j)$ is a `dotted $(K, \epsilon )$-quasigeodesic' for some $K,
\epsilon$, i.e. a quasigeodesic which is not necessarily
connected. Note that by our definition of a quasigeodesic (viz. a
quasi-isometric image of an interval in $\Bbb{R}$) a dotted
quasigeodesic is a quasigeodesic. Since $X$ is $\delta$-hyperbolic,
there exists $C \geq 0$ such that the geodesic $\mu$ lies in a $C$
neighborhood of $\Pi_\lambda (\mu \cap \bigcup_j
Y_j)$. Further, since the latter lies on $B_\lambda$, we conclude that
$B_\lambda$ is $C$-quasiconvex.
$\Box$

\smallskip

\noindent {\bf Note:} Suppose instead of being geodesics $\lambda_i$
are all $(K_Y, \epsilon_Y)$-quasigeodesics in each $Y_i$. Then observe
that the Theorem \ref{mainref} still goes through, and the union
$B_\lambda$ (in this case of {\bf quasigeodesics} ) of $\lambda_i$ is
still a retract of $\bigcup_i Y_i$.

\smallskip

The following lemma gives us a method of finding an admissible path
from $\lambda_j$ to $\lambda_0$, whose length is of the order of $j$.

\begin{lemma}
There exists  $A > 0$, such that 
if $a\in \lambda_j$ for some $j$ then
there exists $b\in{i({\lambda})} = \lambda_0$ such that
${d_X}(a,b)$ $ \leq $ ${A}j$.
\label{connectionlemma}
\end{lemma}
 
\noindent {\bf Proof:}
It suffices to prove that for $j \geq 1$,
there exists $A > 0$ independent of $j$
such that if 
$p \in {\lambda_{j}}$, there exists 
$q \in ({\lambda_{{j-1}}})$
with $d(p,q) \leq A$. A symmetric  argument works for $j \leq -1$.

By construction, $\lambda_{j} = {\Phi_{j}}({\lambda_{j-1}})$.
Since $\phi_{j}$ is a $(K, \epsilon )$ quasi-isometry for all $j$,
there exists $C_1$ such that $\phi_j^{-1} (p)$ 
lies in a $C_1$ neighborhood of 
$\lambda_j$. Hence, there exists $q \in \lambda_j$
such that $d({q},p) \leq 1 + C_1$. Choosing $A = C_1 + 1$ we are through.
 $\Box$

\smallskip

\noindent {\bf Remark:} Note that Lemma \ref{connectionlemma} also
goes through when each $\lambda_i$ is a $(K_Y,
\epsilon_Y)$-quasigeodesic in $Y_i$ for all $i$.

\smallskip

The main theorem of this section follows:

\begin{theorem}
 Let (X,d) be a tree (T) of hyperbolic metric spaces satisfying the
qi-embedded condition.  Suppose in addition that $T$ is either
$\Bbb{R}$ or $[0, \infty )$ with the usual vertex and edge sets $Y_j$,
  for $j$ integers.
If $X$ is hyperbolic, then
${i} : Y_0 \rightarrow X$ extends continuously to ${\hat{{i}}} : \widehat{Y_0}
\rightarrow \widehat{X}$.
\label{main}
\end{theorem}

\noindent {\bf Proof:} 
To prove the existence of a Cannon-Thurston map,
it suffices to show  (from  Lemma \ref{contlemma})
that for all $M\geq{0}$ and ${y_0} \in Y_0$ there exists $N\geq{0}$
such that if a geodesic segment 
$\lambda$ lies outside the $N$-ball around ${y_0}\in{Y_0}$,
$B_{\lambda}$ lies outside the $M$-ball around ${y_0} \in {X}$.

To prove this, we show that if $\lambda$ lies outside the 
$N$-ball around ${y_0}\in Y_0$,
$B_{\lambda}$ lies outside a certain  $M(N)$-ball around $y_0 \in Y_0
\subset {X}$,
 where
$M(N)$ is a proper function from $\Bbb{N}$ into itself.

Since $Y_0$ is properly embedded in $X$ there exists $f(N)$
such that $i( \lambda )$
 lies outside the $f(N)-$ball around $y_0$ 
in $X$
and $f(N)\rightarrow\infty$ as $N\rightarrow\infty$.

Let $p$ be any point on $B_{\lambda}$. Then $p\in Y_j$ for some $j$.
There exists $y\in \lambda =\lambda_0$
such that ${d_X}(y,p)\leq{A}j$ by Lemma \ref{connectionlemma}.
Therefore,
\begin{eqnarray*}
{d_X}({y_0},p) & \geq & {d_X}({y_0},y)-Aj \\
                     & \geq & f(N)-Aj
\end{eqnarray*}

By our choice of metric on $X$,

\begin{center}
${d_X}({y_0},p) \geq j $
\end{center}

Hence
\begin{eqnarray*}
 {d}({x_0},p) & \geq & max({f(N)-Aj}, j) \\
                 & \geq & \frac{f(N)}{A+1}
\end{eqnarray*}

From Theorem \ref{mainref}
there exists $C$ independent of $\lambda$ such that
 $B_{\lambda}$ is a $C$-quasiconvex set containing $\lambda$.
Therefore
any geodesic joining the end-points of $\lambda$ lies in a 
$C$-neighborhood of $B_{\lambda}$.

Hence any geodesic joining  end-points of $\lambda$ lies outside a ball of radius $M(N)$ where
\begin{center}
$M(N) = {\frac{f(N)}{A+1}}-C$
\end{center}
Since $f(N)\rightarrow\infty$ as $N\rightarrow\infty$ so does $M(N)$.
$\Box$

\section{Closed Surface Groups of Bounded Geometry: \\ Theorems of
  Cannon-Thurston and Minsky}

\subsection{Three Manifolds Fibering Over the Circle}

Now, let $M$ be a closed hyperbolic 3-manifold, fibering over the
circle with fiber $F$. Then the universal cover $\tilde{M}$ may be
regarded as a tree $T$ of spaces, where $T = \Bbb{R}$, and where each
vertex and edge set is $\tilde{F}$. Thus $X$ is $\tilde{M}$, which is
quasi-isometric to ${\Bbb{H}}^3$ and $Y$ is $\tilde{F}$ which is
quasi-isometric to ${\Bbb{H}}^2$. Further, $\phi_i$ can all be
identified with $\tilde{\phi}$, where $\phi$ is the pseudo-anosov
monodromy of the fibration. Then as a direct consequence of Theorem
\ref{main}, we obtain the following theorem of Cannon and Thurston,
which motivated much of the present work.

\begin{theorem} {\bf (Cannon-Thurston) \cite{CT}} $M$ be a closed
  hyperbolic 3-manifold, fibering over the 
circle with fiber $F$. Let $i: {\Bbb{H}}^2  \rightarrow
{\Bbb{H}}^3$ denote the inclusion of $\tilde{F}$ into $\tilde{M}$
(where $\tilde{F}$ and $\tilde{M}$ are identified with ${\Bbb{H}}^2$
and
${\Bbb{H}}^3$ respectively. Then $i$ 
extends continuously to a map $\hat{i}: {\Bbb{D}}^2  \rightarrow
{\Bbb{D}}^3$, where ${\Bbb{D}}^2$ and 
${\Bbb{D}}^3$ denote the compactifications of ${\Bbb{H}}^2  $ and $
{\Bbb{H}}^3$ respectively.
\label{ct}
\end{theorem}

\subsection{Closed Surface Kleinian Groups}

In this section we apply Theorem \ref{main} to geometrically tame closed
groups of bounded geometry.
 
The {\it convex core} of a hyperbolic 3-manifold $N$ (without cusps)
is the smallest convex submanifold $C(N) \subset N$ for which inclusion 
is a homotopy equivalence. If an $\epsilon$ neighborhood of
$C(N)$ has finite volume, $N$ is said to be
{\it geometrically finite}. Suppose $N = {\Bbb{H}}^3/ \rho (\pi_1
(S))$ for a closed surface $S$.  We say that an
end of $N$ is {\it geometrically finite} if it has a neighborhood missing
$C(N)$. $N$ is {\it simply degenerate} if it has only one end $E$;
 if  a neighborhood of $E$ is
homeomorphic to $S{\times}{\Bbb{R}}$; and if, in addition,
 there is a sequence of
pleated surfaces 
homotopic in this neighborhood to the inclusion of $S$, and exiting every
compact set (This last condition is automatic by Bonahon
\cite{bonahon-bouts}. $N$ is called {\it doubly degenerate} if it has
two ends, both of which are 
 simply degenerate. For a more detailed discussion of pleated surfaces 
and geometrically tame ends, see \cite{Thurstonnotes} or \cite{minsky-top}.

Let ${inj}_N(x)$ denote the injectivity radius at $x\in{N}$. 
$N$ is said to have {\bf bounded geometry} if  there exists
$\epsilon_0 > 0$ such that 
${inj}_N{(x)} > \epsilon_0$ for all $x\in{N}$. 
In order to apply Theorem \ref{main}
we need some preliminary Lemmas.

Let $E$ be a simply degenerate end of $N$. Then $E$ is homeomorphic to
$S{\times}[0,{\infty})$ for some closed surface $S$ of genus greater than one.

\begin{lemma}
\cite{Thurstonnotes}
There exists $D_1 > 0$ such that for all $x\in{N}$, there exists a pleated 
surface $g : (S,{\sigma}) \rightarrow N$ with $g(S){\cap}{B_{D_1}}(x) \neq \emptyset$.
\label{closepleated}
\end{lemma}

The following Lemma follows easily from the fact that ${inj}_N{(x)} > \epsilon_0$:

\begin{lemma} \cite{bonahon-bouts},\cite{Thurstonnotes}
There exists $D_2 > 0$ such that if $g : (S,{\sigma}) \rightarrow N$ is a
pleated surface, then $dia(g(S)) < D_2$. 
\label{diameter}
\end{lemma}

The following Lemma due to Minsky \cite{minsky-top} follows from compactness
of pleated surfaces.

\begin{lemma}
\cite{minsky-top}
Fix $S$ and $\epsilon > 0$. Given $a > 0$ there exists $b > 0$ such that if
$g : (S,{\sigma})\rightarrow{N}$
and $h : (S,{\rho})\rightarrow{N}$ are homotopic pleated surfaces which
are isomorphisms on $\pi_1$ and ${inj}_N{(x)} > \epsilon$ for all $x \in N$,
then\\
\begin{center}
${d_N}(g(S),h(S)) \leq a \Rightarrow {d_{Teich}}({\sigma},{\rho}) \leq b$,
\end{center}
    where $d_{Teich}$ denotes Teichmuller distance.
\label{pleatedcpt}
\end{lemma}

{\bf Definition:} {\it The {\bf universal curve} over $X{\subset}Teich(S)$
is a bundle whose fiber over $x\in{X}$ is $x$ itself.} \cite{ctm-locconn}

\begin{lemma}
There exist $K, \epsilon$ and  a homeomorphism $h$ from $E$ to the
universal curve $S_\gamma$
over a Lipschitz path $\gamma$ in Teichmuller space, such that
$\tilde{h}$ from $\tilde{E}$ to the universal cover of $S_\gamma$  
is a $(K,{\epsilon})$-quasi-isometry.
\label{lipcurve}
\end{lemma}

{\bf Proof:}
We can assume that $S{\times}\{{0}\}$ is mapped to a pleated surface
$S_0$ $\subset N$ under the homeomorphism from $S{\times}[0,{\infty})$
to $E$. We shall construct inductively a sequence of `equispaced' pleated
surfaces ${S_i}\subset{E}$ exiting the end. Assume that ${S_0},{\cdots},{S_n}$
have been constructed such that:
\begin{enumerate}
\item If $E_i$ be the non-compact component of $E{\setminus}{S_i}$, then
$S_{i+1} \subset E_i$.
\item Hausdorff distance between $S_i$ and $S_{i+1}$ is bounded above by
$3({D_1}+{D_2})$.
\item ${d_N}({S_i},{S_{i+1}}) \geq D_1 + D_2$.
\item From Lemma \ref{pleatedcpt} and condition (2)
above there exists $D_3$ depending on $D_1$, $D_2$ and $S$ such that
$d_{Teich}({S_i},{S_{i+1}}) \leq D_3$
\end{enumerate}

Next choose $x \in E_n$, such that ${d_N}(x,{S_n}) = 2({D_1}+{D_2})$. Then
by Lemma \ref{closepleated}, there exists a pleated surface
$g : (S,{\tau}) \rightarrow N$ such that 
${d_N}(x,{g(S)}) \leq D_1$. Let ${S_{n+1}} = g(S)$. Then by the triangle
inequality and Lemma \ref{diameter}, if $p\in{S_n}$ and $q\in{S_{n+1}}$,
\begin{center}
${D_1} + {D_2} \leq {d_N}(p,q) \leq 3({D_1} + {D_2})$.
\end{center}

This allows us to continue inductively. The Lemma follows. $\Box$ 

\medskip

Note that in the above Lemma, pleated surfaces are not assumed to be embedded.
This is because   immersed pleated surfaces with a uniform lower
bound on  injectivity
radius are uniformly quasi-isometric to the corresponding
Riemann surfaces.

 Since there are exactly one or two
 ends $E_i$, we have thus shown:

\begin{lemma}
The hyperbolic metric space $\widetilde{C(N)}$ is quasi-isometric to a tree
(T) of hyperbolic metric spaces satisfying the qi-embedded condition,
where $T$ is either $[0, \infty )$ (simply degenerate) or $\Bbb{R}$
  (doubly degenerate). 
\label{hyptree}
\end{lemma}

Applying Theorem \ref{main}, we obtain the following theorem of Minsky:

\begin{theorem} {\bf Minsky \cite{minsky-jams}}
 Let $\Gamma = \rho (\pi_1 (S))$ be a closed surface Kleinian 
group, such that ${{\Bbb{H}}^3}/{\Gamma} = M$ has injectivity radius 
uniformly bounded below by some $\epsilon > 0$. Then there exists a continuous
map from the Gromov boundary of $\Gamma$ (regarded as an abstract group)
to the limit set of $\Gamma$ in ${\Bbb{S}}^2_{\infty}$.
\label{minsky}
\end{theorem}

Since a continuous image of a compact locally connected set is locally 
connected \cite{hock-young},
and since the limit set of $\Gamma$ is a continuous image of the
circle, by Theorem \ref{main}, we have: \\

\begin{cor}
Suppose $\Gamma$ is a closed surface Kleinian group, such that
 $N = {{\Bbb{H}}^3}/{\Gamma}$ has bounded geometry, i.e. 
${inj}_N{(x)} > \epsilon_0$ for all $x \in N$. Then the limit set of 
$\Gamma$ is locally connected.
\label{locconn}
\end{cor}

\section{Punctured
 Surface Groups of Bounded Geometry:  A Theorem of Bowditch}

\subsection{Outline of Proof}
In 
\cite{bowditch-ct} \cite{bowditch-stacks},  Bowditch
proved the existence of Cannon-Thurston maps for punctured surface groups
 of bounded geometry using some of the ideas from  \cite{mitra-trees}. 
We give below a different proof of the result, which is in some
 ways simpler. First, a sketch. 

Let $N^h$ be a bounded geometry 3-manifold corresponding to a
representation of the fundamental group of a punctured surface $S^h$.
Excise
the cusps (if any) of $N^h$ leaving us a manifold that has one or two ends. Let
$N$ denote $N^h$ minus cusps. Then $N$ is quasi-isometric to the
universal curve over a Lipschitz path in Teichmuller space from which
cusps have been removed.  This path is semi-infinite
or bi-infinite 
according as $N$ is one-ended or two-ended. Fix a reference finite
 volume hyperbolic surface $S^h$. Let $S$ denote $S^h$ minus 
cusps. Then $\widetilde{S}$ is quasi-isometric to the Cayley graph of 
$\pi_1{(S)}$ which is (Gromov) hyperbolic. We fix a base surface in $N$
and identify it with $S$. Now look at
$\widetilde{S}\subset\widetilde{N}$. 
 Let $\lambda = [a,b]$ be a geodesic segment in $\widetilde{S}$.
We `flow' $\lambda$ out the end(s) of $\widetilde{N}$ to generate
the {\bf hyperbolic ladder-like set} $B_{\lambda}$ as in the proof of
Theorem \ref{mainref}. 

 Theorem \ref{mainref}
 ensures that there is a retraction from  $\widetilde{N}$ to
$B_\lambda$ which does not increase distances much. From this it
 follows that $B_\lambda$ is 
quasi-isometrically embedded in $\widetilde{N}$.   Recall that for the
 construction of 
$B_\lambda$, we only require the hyperbolicity of $\tilde{S}$ and not
that of $\tilde{N}$. 

Now if $\lambda$ lies outside a large ball about a fixed reference
point $p$ in  
$\widetilde{S}$, then so does $B_\lambda$ in $\widetilde{N}$.  Since
$B_\lambda$ is q.i. embedded in $\widetilde{N}$, there exists an 
ambient 
$\widetilde{N}$-quasigeodesic $\mu$ lying in a bounded neighborhood of 
$B_\lambda$ joining the end-points of $\lambda$. If $S^h$ had no
cusps, 
as in the case for closed surfaces, we could immediately
conclude that for any geodesic segment $\lambda$ in $\widetilde{S^h}$ lying outside
large balls around $p$, there is a quasigeodesic in $\widetilde{N^h}$ joining its
endpoints and lying outside a large ball around $p$  in
$\widetilde{N^h}$. 
However, since
 $S^h$ has cusps, $\widetilde{S^h}$ and $\widetilde{S}$ are different. So a 
little more work is necessary. Suppose as before that $\lambda_0$ is a geodesic in
$\widetilde{S^h}$ lying outside a large ball around $p$. For ease of exposition we assume
that the end-points of $\lambda_0$ lie outside cusps.
 Let $\lambda \subset \widetilde{S}$ be the geodesic in (the path-metric on)
$\widetilde{S}$ joining the same pair of points. Then off horodisks,
$\lambda_0$ and $\lambda$ track each other. Construct 
$B_\lambda$ as before, and let
$\mu$ be an ambient quasigeodesic 
in $\widetilde{N}$ lying in a bounded neighborhood of 
$B_\lambda$ joining the end-points of $\lambda$. Then, modulo horoballs in 
$\widetilde{N^h}$, $\mu$ lies outside a large ball around $p$. Let $\mu_0$ be the
hyperbolic geodesic joining the end points of $\mu$. Off horoballs, $\mu$ and
$\mu_0$ track each other. Hence, off horoballs, $\mu_0$ lies outside large balls
about $p$. The points at which $\mu_0$ enters and leaves
 a particular horoball therefore lie 
 outside large balls 
about $p$. But then the hyperbolic segment joining them must do the
same. This shows 
that $\mu_0$ must itself lie outside large balls around $p$. As before
we conclude 
that there exists  a continuous
extension of the inclusion of  $\widetilde{S^h}$ into $\widetilde{N^h}$ to
the boundary. The remaining part of this section fleshes out this
argument.

\subsection{Relative Hyperbolicity}

We shall  be requiring certain properties of hyperbolic spaces
minus horoballs. These were studied by Farb \cite{farb-relhyp} under
the garb of `electric geometry'. We
combine Farb's results with a 
version that is a (slight variant of) theorem due to McMullen
(Theorem 8.1 of \cite{ctm-locconn}).

{\bf Definition:} A path $\gamma : I \rightarrow Y$ to a path metric space $Y$ is an ambient
K-quasigeodesic if we have
\begin{center}
$L({\beta}) \leq K L(A) + K$
\end{center}
for any subsegment $\beta = \gamma |[a,b]$ and any path $A : [a,b] \rightarrow Y$ with the
same endpoints.

The following definitions are adapted from \cite{farb-relhyp}

{\bf Definition:} Let $M$ be a convex hyperbolic manifold.
Let $Y$ be the universal cover of  $M$ minus cusps
and $X = \widetilde{M}$.  $\gamma$ is said to be a $K$-quasigeodesic in $X$ 
{ \bf without backtracking } if \\ 
$\bullet$ $\gamma$ is  a
  $K$-quasigeodesic in $X$ \\
$\bullet$  $\gamma$
does not return to any
  horoball $\bf{H}$ after leaving it. \\

{\bf Definition:}  $\gamma$ is said to be an ambient
$K$-quasigeodesic in $Y$
{ \bf without backtracking } if \\ 
$\bullet$ $\gamma$ is an ambient $K$-quasigeodesic in $Y$ \\
$\bullet$ $\gamma$ is obtained from a $K$-quasigeodesic without
backtracking in $X$ by
replacing each  maximal subsegment with end-points on a horosphere by
a quasigeodesic lying on the surface of the horosphere. \\

Note that in the above definition, we allow the behavior to be quite
arbitrary on horospheres (since Euclidean quasigeodesics may be quite
wild); however, we do not allow wild behavior off horoballs.

\smallskip

\noindent {\bf Remark:} Our definition of {\it ambient quasigeodesic
  without backtracking} does not allow a path to follow a horosphere
  for a long distance without entering it. This is a point where the
  definition differs from the usual definition of an ambient
  quasigeodesic. 

\smallskip

$B_R (Z)$ will denote the $R$-neighborhood of the set $Z$. \\
Let $\cal{H}$ be a locally finite collection of horoballs in a convex
subset $X$ of ${\Bbb{H}}^n$ 
(where the intersection of a horoball, which meets $\partial X$ in a point, 
 with $X$ is
called a horoball in $X$). 
The following theorem is due to McMullen \cite{ctm-locconn}.

\begin{theorem} \cite{ctm-locconn} 
Let $\gamma: I \rightarrow X \setminus \bigcup \cal{H} = Y$ be an ambient
$K$-quasigeodesic for a convex subset $X$ of ${\Bbb{H}}^n$ and let
$\mathcal{H}$  denote a
collection of horoballs.
Let $\eta$ be the hyperbolic geodesic with the same endpoints as
$\gamma$. Let $\cal{H}({\eta})$  
be the union of all the horoballs in $\mathcal{H}$ meeting $\eta$. Then
$\eta\cup\mathcal{H}{({\eta})}$ is (uniformly) quasiconvex and $\gamma
(I) \subset  
B_R (\eta \cup \cal{H} ({\eta}))$, where $R$ depends only on
$K$. 

\label{ctm}
\end{theorem}

Theorem \ref{ctm} is similar in flavor to certain theorems
 about relative hyperbolicity
{\it a la} Gromov \cite{gromov-hypgps}, Farb \cite{farb-relhyp} and Bowditch
\cite{bowditch-relhyp}. We give below a related theorem that is
derived from Farb's
`Bounded Horosphere Penetration' property. 

Let $\gamma_1 = \overline{pq}$ be a hyperbolic $K$-quasigeodesic
without backtracking starting from a horoball $\bf{H_1}$ and ending
within (or on) a {\em different}
horoball ${\bf{H_2}}$. Let $\gamma = [a,b]$ be the hyperbolic geodesic
minimizing distance between $\bf{H_1}$ and $\bf{H_2}$. Following
\cite{farb-relhyp} we put the zero metric on the horoballs that
$\gamma$ penetrates. The resultant pseudo-metric is called the
electric metric. Let $\widehat{\gamma}$ and  $\widehat{\gamma_1}$
 denote the paths represented by $\gamma$ and $\gamma_1$
respectively in this pseudometric. It is shown in \cite{farb-relhyp}
that $\gamma$, $\widehat{\gamma}$ and  $\widehat{\gamma_1}$ have
similar intersection patterns with horoballs, i.e. 
there exists $C_0$ such that \\
$\bullet 1$ If only one of $\gamma$ and $\widehat{\gamma_1}$
 penetrates a horoball
  $\bf{H}$, then it can do so for a distance $ \leq C_0$. \\
$\bullet 2$ If both $\widehat{\gamma_1}$ and $\gamma$ enter (or leave) a horoball
  $\bf{H}$ then their entry (or exit) points are at a distance of at
  most $C_0$ from each other. [Here by `entry' (resp. `exit') point of a
  quasigeodesic we mean
  a point at which the path switches from being in the complement of
  or `outside'
 (resp. in the interior of or `inside') a closed horoball to being
  inside (resp. outside) it].\\
The point to observe here is that quasigeodesics without backtracking
in our definition gives rise to quasigeodesics without backtracking in
 Farb's sense. Since this is true for arbitrary $\gamma_1$ we give
 below a slight strengthening of this fact. Further, by our
 construction of ambient quasigeodesics without backtracking, we might
 just as well consider ambient quasigeodesics without backtracking in
 place of quasigeodesics.

\begin{theorem}  \cite{farb-relhyp}
Given $C > 0$, there exists $C_0$ such that if \\
$\bullet 1$ either two quasigeodesics without backtracking 
$\gamma_1 , \gamma_2$ in $X$, OR\\
$\bullet 2$ two ambient
quasigeodesics without backtracking  $\gamma_1 , \gamma_2$ in $Y$, OR \\
$\bullet 3$ $\gamma_1$ - an ambient
quasigeodesic without backtracking in $Y$ and $\gamma_2$ -
a quasigeodesic without
backtracking in $X$,\\

\smallskip

\noindent start and end \\

\smallskip

\noindent $\bullet 1$ either  on (or within) the same horoball OR\\
$\bullet 2$  a distance $C$ from each other \\
then they have similar intersection patterns with horoballs (except
possibly the first and last ones), i.e.
there exists $C_0$ such that \\
$\bullet 1$ If only $\gamma_1$ 
penetrates or travels along the boundary of a horoball
  $\bf{H}$, then it can do so for a distance $ \leq C_0$. \\
$\bullet 2$ If both $\gamma_1$ and $\gamma_2$ enter (or leave) a horoball
  $\bf{H}$ then their entry (or exit) points are at a distance of at
  most $C_0$ from each other. 

\label{farb}
\end{theorem}

\subsection{Horo-ambient Quasigeodesics}

A special kind of 
quasigeodesic without back-tracking will be necessary.
We start with a
hyperbolic geodesic $\lambda^h$ in ${\widetilde{S}}^h$. Fix a
neighborhood of the 
cusps lifting to an equivariant family of horoballs in the universal
cover ${\Bbb{H}}^2 = \widetilde{S^h}$. 
 Since $\lambda^h$ is a hyperbolic
geodesic in  $\widetilde{S^h}$ there are unique entry and exit
points for each horoball that $\lambda^h$ meets and hence unique
Euclidean geodesics joining them on the corresponding
horosphere. Replacing the segments of $\lambda^h$ lying inside
$Z$-horoballs by the corresponding Euclidean geodesics, we obtain an
ambient quasigeodesic $\lambda$ in $\widetilde{M_0}$ as a consequence
of
 Theorem
\ref{ctm} (See Corollary \ref{ctm-cor} and figure  below):

\medskip

\begin{center}
\includegraphics{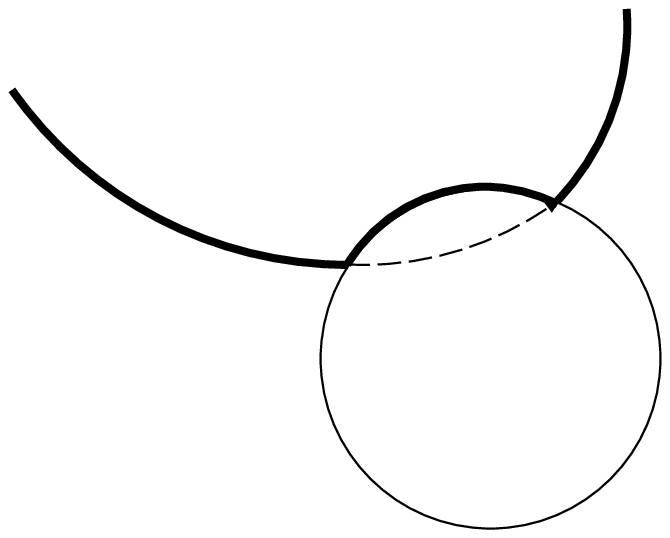}

\smallskip

\underline{\it Horo-ambient quasigeodesics}

\end{center}

\smallskip

Ambient quasigeodesics obtained by this kind of a construction will be
termed {\bf horo-ambient quasigeodesics} to distinguish them from {\em
  electro-ambient quasigeodesics} defined earlier. 

The following Corollary of Theorem \ref{ctm} justifies the
terminology. 

\begin{cor} 
There exists $K, \epsilon$, such that any horo-ambient quasigeodesic
in $Y$ is indeed a $K, \epsilon$
ambient quasigeodesic in $Y$.
\label{ctm-cor}
\end{cor}

\noindent {\bf Proof:} Let $[a,b]_h$ be a hyperbolic geodesic in $X$
joining $a, b$, where $a, b$ lie outside horoballs. Let $[a,b]_{ha}$ be
the horo-ambient quasigeodesic in $Y$ joining $a, b$. Let $[a,b]_a$ be
the ambient geodesic in $Y$ joining $a, b$. Then by Theorem
\ref{ctm},there exists $R \geq
0$, such that  $[a,b]_a$ lies in an $R$-neighborhood of $[a,b]_h \cup
\mathcal{H} (\eta )$. Project $[a,b]_a$ onto $[a,b]_h \cup
\mathcal{H} (\eta )$ using the nearest point projection in
$X$. Removing the back-tracking induced, we get some
$K, \epsilon$ depending on $R$ such that the image is 
an ambient quasigeodesic without backtracking in our sense. Clearly,
such an ambient quasigeodesic coincides with $[a,b]_{ha}$ off
horoballs $\mathcal{H} (\eta )$. Since the interpolating segments in
$[a,b]_{ha}$ on
any horoball in $\mathcal{H} (\eta )$ are Euclidean {\em geodesics},
the result follows. $\Box$

\subsection{Trees of Spaces}

We want to first show that the universal cover of $N^h$ minus $Z$-cusps is
quasi-isometric to a tree of hyperbolic metric spaces.

Let $E^h$ be a simply degenerate end of $N^h$. Then $E^h$ is homeomorphic to
$S^h{\times}[0,{\infty})$ for some  surface $S^h$ of negative Euler
  characteristic.
Cutting off a neighborhood of the cusps of $S^h$ we get a surface with boundary
denoted as $S$. Let $E$ denote $E^h$ minus a neighborhood of the
 $Z$-cusps. We assume that each $Z$-cusp has the standard form coming from
a quotient of a horoball in ${\Bbb{H}}^3$ by $Z$. Also, we shall take
our pleated 
surfaces to be such that the pair $(S,cusps)$ is mapped to the pair
$(E,cusps)$ for each 
pleated $S^h$. We shall now show that each $\widetilde{E}$ is 
quasi-isometric to a ray of hyperbolic metric spaces satisfying the q-i
embedded condition.  Each edge and vertex space will be a copy of 
$\widetilde{S}$ and the edge to vertex space inclusions shall be
quasi-isometries. Note that each $\widetilde{S}$ can be thought of as 
a copy of ${\Bbb{H}}^2$ minus an equivariant family of horodisks.
The following Lemmas are  generalizations to punctured surfaces of
Lemmas \ref{closepleated}, \ref{diameter}, \ref{pleatedcpt}. 

\begin{lemma}
\cite{Thurstonnotes}
There exists $D_1 > 0$ such that for all $x \in {E}$,
 there exists a pleated 
surface $g : (S^h,{\sigma}) \rightarrow E^h$ with 
$g(S){\cap}{B_{D_1}}(x) \neq \emptyset$. Also $g$ maps $(S, cusps)$
to $(E, cusps)$.
\label{closepleated-punct}
\end{lemma}

\begin{lemma} \cite{bonahon-bouts},\cite{Thurstonnotes}
There exists $D_2 > 0$ such that if $g : (S^h,{\sigma}) \rightarrow N^h$ is a
pleated surface, then the diameter of the image of $S$ is bounded,
i.e.
$dia(g(S)) < D_2$. 
\label{diameter-punct}
\end{lemma}

The following Lemma due to Thurston (Theorems 9.2 and 9.6.1 of
\cite{Thurstonnotes} ) and Minsky \cite{minsky-top} follows from compactness
of pleated surfaces.

\begin{lemma}
\cite{minsky-top}
Fix $S^h$ and $\epsilon > 0$. Given $a > 0$ there exists $b > 0$ such that if
$g : (S^h,{\sigma})\rightarrow{E^h}$
and $h : (S^h,{\rho})\rightarrow{E^h}$ are homotopic pleated surfaces which
are isomorphisms on $\pi_1$ and $E^h$ is of bounded geometry,
then\\
\begin{center}
${d_E}(g(S),h(S)) 
\leq a \Rightarrow {d_{Teich}}({\sigma},{\rho}) \leq b$,
\end{center}
    where $d_{Teich}$ denotes Teichmuller distance.
\label{pleatedcpt-punct}
\end{lemma}

In \cite{minsky-top} a specialization of this statement is proven for
closed surfaces. However, the main ingredient, a Theorem due to
Thurston is stated and proven in \cite{Thurstonnotes} (Theorems 9.2
and 9.6.1 - {\it 'algebraic limit is geometric limit'}) for finite area
surfaces. The arguments given by Minsky to prove the above Lemma from
Thurston's Theorems (Lemma 4.5, Corollary 4.6 and Lemma 4.7 of
\cite{minsky-top}) go through with very little change for surfaces
of finite area. 

\smallskip

{\bf Construction of equispaced pleated surfaces exiting the end}\\
We next construct a sequence of {\bf equispaced} pleated surfaces
$S^h(i)\subset E^h$ exiting the end as before. 
 Assume that ${S^h(0)},{\cdots},{S^h(n)}$
have been constructed such that:\\
\begin{enumerate}
\item  $(S(i),cusps)$ is mapped to $(E, cusps)$ \\
\item If $E(i)$ be the component of $E{\setminus}{S(i)}$ for which
$E(i)$ is non-compact, then
$S(i+1) \subset E(i)$. \\
\item Hausdorff distance between $S(i)$ and $S(i+1)$ is bounded above by
$3({D_1}+{D_2})$. \\
\item ${d_E}({S(i)},{S(i+1)}) \geq D_1 + D_2$. \\
\item From Lemma \ref{pleatedcpt-punct} and condition (3)
above there exists $D_3$ depending on $D_1$, $D_2$ and $S$ such that
$d_{Teich}({S(i)},{S({i+1})}) \leq D_3$ \\
\end{enumerate}

Next choose $x \in E(n)$, such that ${d_E}(x,{S_n}) = 2({D_1}+{D_2})$. 
Then
by Lemma \ref{closepleated-punct}, there exists a pleated surface
$g : (S^h,{\tau}) \rightarrow E^h$ such that 
${d_E}(x,{g(S)}) \leq D_1$. Let ${S^h(n+1)} = g(S^h)$. Then by the triangle
inequality and Lemma \ref{diameter-punct}, if $p\in{S(n)}$ and 
$q\in{S(n+1)}$,
\begin{center}
${D_1} + {D_2} \leq {d_E}(p,q) \leq 3({D_1} + {D_2})$.
\end{center}

This allows us to continue inductively. 
$S(i)$ corresponds to a point $x_i$ of $Teich(S)$. Joining the $x_i$'s in
order, one gets a Lipschitz path in $Teich(S)$. 

\smallskip

{\bf Definition:} A sequence of pleated surfaces satisfying conditions
(1-5) above will be called an {\bf  equispaced sequence of pleated surfaces}.
The corresponding sequence of $S(i) \subset E$ will be called 
an {\bf equispaced sequence of  truncated pleated surfaces}

Each $S(i)$ being compact (with or without boundary),
$\widetilde{S(i)}$ is a hyperbolic metric space. 
We can think of
the universal cover
$\widetilde{E}$ of $E$ as being quasi-isometric to a  ray $T$
of hyperbolic metric spaces by setting $T = [0,{\infty})$, with vertex
set $\mathcal{V} = \{$ $ n : n \in {\Bbb{N}}{\cup}$ $\{$ $ 0 $ $\}$ $ \}$, 
edge set 
$\mathcal{E} =$ $ \{$ $ [{n-1},n]: n \in {\Bbb{N}} $ $ \}$,
${X_n} = {\widetilde{S(n)}} = X_{[n-1,n]}$.  
Further, by Lemma \ref{pleatedcpt-punct}
this tree of hyperbolic metric spaces satisfies the quasi-isometrically
embedded condition. We have thus shown

\begin{lemma}
If $E^h$ be a simply degenerate  end of a hyperbolic 3 manifold $N^h$
with bounded geometry, then there is a sequence of equispaced 
pleated surfaces exiting $E^h$ and hence a sequence of truncated
pleated surfaces exiting $\widetilde{E}$. Further, $\widetilde{E}$ is
quasi-isometric to a ray of hyperbolic metric spaces satisfying the 
q.i. embedded condition.
\label{equispaced-punct}
\end{lemma}

\subsection{Construction of $B_\lambda$: Modifications for Punctured Surfaces}

Now, let $N^h$ be a bounded geometry 3-manifold corresponding to a
representation of the fundamental group of a punctured surface $S^h$.
Excise
the cusps (if any) of $N^h$ leaving us a manifold that has one or two ends. Let
$N$ denote $N^h$ minus cusps. 
Fix a reference finite
 volume hyperbolic surface $S^h$. Let $S$ denote $S^h$ minus 
cusps. Then $\widetilde{S}$ is quasi-isometric to the Cayley graph of 
$\pi_1{(S)}$ which is (Gromov) hyperbolic, in fact quasi-isometric to
a tree.  We fix a base surface in $N$
and identify it with $S$. Now look at
$\widetilde{S}\subset\widetilde{N}$. 

Then, by Lemma \ref{equispaced-punct} $\widetilde{N}$ is
quasi-isometric to a tree $T$ of hyperbolic metric spaces. Each of the
vertex and edge spaces
 is a copy of $\widetilde{S}$. Also, the map $\phi_i$
induced from $\widetilde{S} \times \{ i \}$ to
 $\widetilde{S} \times \{ i+1 \}$ is a $(K,\epsilon )$-quasi-isometry
for all $i$ as are their quasi-isometric inverses.
Further, $T$ is a semi-infinite
or bi-infinite interval in $\Bbb{R}$
according as $N$ is one-ended or two-ended. So far, this is exactly
like the case for closed surfaces. But here, we can  assume 
in addition that each $\phi_i$ is the restriction of a map
$\phi_i^h$ from  ${\widetilde{S}}^h \times \{ i \}$
to  ${\widetilde{S}}^h \times \{ i+1 \}$ which preserves horodisks.
Let $\Phi_i^h$ denote the induced map on hyperbolic geodesics and let
$\Phi_i$ denote the induced map on horo-ambient geodesics.

Let $\lambda^h$ be a hyperbolic geodesic segment in
${\widetilde{S}}^h$. Let $\lambda$ be the horo-ambient quasigeodesic
in $\widetilde{S}$ joining the end-points of $\lambda^h$.

Starting with a horo-ambient quasigeodesic $\lambda \subset
\widetilde{S}$, we can now proceed as in the proof of Theorem
\ref{mainref}
 to construct  the
 {\bf hyperbolic ladder-like set} 
 $B_{\lambda}$.

There is only one difference: {\bf Each $\lambda_i$ in this situation
 is a horo-ambient quasigeodesic, and not necessarily a hyperbolic
 geodesic.} 
Thus, we set $\lambda = \lambda_0$ to be some horo-ambient
 quasigeodesic in $\widetilde{S} = \widetilde{S} \times \{ 0
 \}$. Next, (for $i \geq 0$),
inductively, set $\lambda_{i+1}$ to be the (unique)
 horo-ambient quasigeodesic in $\widetilde{S} \times \{ i+1 \}$
 joining the end-points of $\phi_i (\lambda_i )$. That is to say,
 $\lambda_{i+1} = \Phi_i (\lambda_i )$. 
Similarly, for $i \leq 0$. 

Note that Corollary \ref{ctm-cor} ensures
 that there exist $K_0, \epsilon_0$ such that
each $\lambda_i$ is a $(K_0, \epsilon_0 )$-quasigeodesic in $\widetilde{S}
 \times \{ i \}$.

Then  Theorem \ref{mainref} (or more precisely, the Note following it)
 ensures that there is a retraction from  $\widetilde{N}$ to
$B_\lambda$ which does not increase distances much. From this it
 follows that $B_\lambda$ is 
quasi-isometrically embedded in $\widetilde{N}$.   Recall that for the
 construction of 
$B_\lambda$, we only require the hyperbolicity of $\tilde{S}$ and not
that of $\tilde{N}$. 

As before, (by projecting a geodesic in $\widetilde{N}$ onto
$B_\lambda$)
we obtain an ambient quasigeodesic contained in $B_\lambda$ joining
the end-points $a, b$  of $\lambda$.

Let \\
$\bullet$  $\beta^h$ $ = $ geodesic in $\widetilde{N^h}$ joining $a,
b$ \\
$\bullet$   $\beta_{amb}^0$  $ = $ horo-ambient quasigeodesic 
 in $\widetilde{N}$ obtained from  $\beta^h$ by replacement of
 hyperbolic by `Euclidean' geodesic segments for horoballs in
$\widetilde{N^h}$ \\
$\bullet$ $\beta_{amb} = \Pi_\lambda ( \beta_{amb}^0  ) $\\

\subsection{Quasigeodesic Rays}

Let
$\lambda_i^c$ denote the union of the segments of $\lambda_i$ which lie
along horocycles and let $\lambda_i^b = \lambda_i - \lambda_i^c$. Let

\begin{center}

$B^c_\lambda = \cup_i \lambda_i^c $ \\
$B^b_\lambda = \cup_i \lambda_i^b $\\

\end{center}

We want  to show that for all $x \in B^b_\lambda$ there exists a
$C$-quasigeodesic $r_x : \{ 0 \} \cup \Bbb{N} \rightarrow B^b_\lambda$
such that $x \in r_x ( \{ 0 \} \cup \Bbb{N} )$ and $r_x (i) \in
\lambda_i^b$. Suppose $x \in \lambda^b_k \subset B^b_\lambda$. We define $r_x$
by starting with $r_x (k) = x$ and construct $r_x (k-i)$ and $r_x
(k+i)$ inductively (of course $(k-i)$ stops at $0$ for $T$ a
semi-infinite ray but $(k+i)$ goes on
to infinity). For the sake of concreteness, we prove the existence of
such a $r_x (k+1)$. The same argument applies to $(k-1)$ and
inductively for the rest.

\begin{lemma}
There exists $C > 0$ such that if $r_x (k) = x \in \lambda_k^b$ then there
exists $x^{\prime} \in \lambda^b_{k+1}$ such that $d(x, x^{\prime} ) \leq
C$. We denote $r_x (k+1) = x^{\prime}$.
\label{qgray1}
\end{lemma}

{\bf Proof:} Let $[a,b]$ be the maximal connected component of
$\lambda_k^b$  on which $x$ lies. Then there exist two horospheres
$\bf{H}_1$ and $\bf{H}_2$ such that  $a \in \bf{H}_1$ (or is the
initial point of $\lambda_k$ ) and  $b \in \bf{H}_2$ (or is the terminal
point of $\lambda_k$ ).  Since $\phi_{k}$ preserves
horocycles, $\phi_{k} (a)$ lies on a horocycle (or is the initial
point of $\lambda_{k+1}$ ) as does $\phi_{k} (b)$ (or is the terminal
point of $\lambda_{k+1}$). Further, the image of $[a,b]$ under
$\phi_k$
is a 
quasigeodesic in $\widetilde{S} \times \{ k+1 \}$
{\bf which we now denote as $\phi_{k} ([a,b])$}.
Recall that $\Phi_{k} ([a,b])$ is the horo-ambient 
geodesic in $\widetilde{S} \times \{ k+1 \}$ 
joining $\phi_{k} (a)$ and $\phi_{k} (b)$. 

Therefore
by (Gromov) hyperbolicity of $\widetilde{S}$,
$\Phi_{k} ([a,b])$ lies in a bounded neighborhood of  $\phi_{k}
([a,b])$ (which in turn lies at a bounded distance from
$\Phi^h_{k}([a,b])$) and hence by Theorem \ref{farb} there exists an
upper bound on 
how much $\Phi^h_{k} ([a,b])$ can penetrate horoballs, i.e. there
exists $C_1 > 0$  such that for all $z \in \Phi^h_{k} ([a,b])$,
there exists $z^{\prime} \in \Phi^h_{k} ([a,b])$ lying outside
horoballs with $d(z, z^{\prime}) \leq C_1$. Further, since
$\phi_{k}^h$ is a quasi-isometry     there exists $C_2 > 0$ such that
$d( \phi_{k} (x), \Phi^h_{k} ([a,b])) \leq C_2$. Hence there
exists $x^{\prime} \in  \Phi^h_{k} ([a,b])$ such that \\
$\noindent \bullet$ $d( \phi_{k} (x) , x^{\prime}) \leq C_1 + C_2$
\\
$\noindent \bullet$ $x^{\prime}$ lies outside horoballs. \\

Again,   $\Phi_{k} ([a,b])$
 lies at a uniformly bounded distance $\leq C_3$ from $\lambda_{k+1}$ and
 so, if $c, d \in \lambda_{k+1}$ such that $d(a,c) \leq C_3$ and $d(b,d)
 \leq C_3$ then the segment $[c,d]$ can penetrate only a bounded
 distance into any horoball. Hence there exists $C_4 > 0$ and
$x^{\prime \prime} \in
 [c,d] \subset \lambda_{k+1}$ such that \\
$\noindent \bullet$ $ d( x^{\prime}, x^{\prime \prime} ) \leq C_4$ \\
$\noindent \bullet$ $ x^{\prime \prime}$ lies outside horoballs. \\

Hence, $d(\phi_{k} (x),  x^{\prime \prime}) \leq C_1 + C_2 +
C_4$. Since $d(x, \phi_{k+1} (x)) = 1$, we have, by choosing
$r_{k+1}(x) =  x^{\prime \prime}$,

\begin{center}

$d( r_k(x), r_{k+1}(x)) \leq 1 + C_1 + C_2 + C_4$.

\end{center}

Choosing $C =  1 + C_1 + C_2 + C_4$, we are through. $\Box$

\medskip

Using Lemma \ref{qgray1} repeatedly (inductively replacing $x$ with
$r_x (k + i)$ we obtain the values of $r_x (i)$ for $i \geq k$. By an
exactly similar symmetric argument, we get $r_x (k-1)$ and proceed
down to $r_x (0)$. Now for any $i$, $z \in {\widetilde{S}} \times \{ i
\}$ and
 $y \in {\widetilde{S}} \times \{ i
+ 1 \}$, $d(z, y) \geq 1$. Hence, for any 
 $z \in {\widetilde{S}} \times \{ i
\}$ and
 $y \in {\widetilde{S}} \times \{ j
\}$, $d(z, y) \geq |i-j|$. This gives

\begin{cor} 
There exist $K, \epsilon > 0$ such that for all $x \in \lambda_k^b \subset
B^b_\lambda$ there exists a $(K, \epsilon )$ quasigeodesic ray $r_x$
such that $r_x (k) = x$ and $r_x (i) \in \lambda_i^b$ for all $i$.
\label{qgray2}
\end{cor}
 
\smallskip

To fix and recall notation:  \\
$\bullet$ $\lambda^h$ $= $ hyperbolic geodesic in ${\widetilde{S}}^h$
joining $a, b$
\\
$\bullet$ $\lambda $ $ = $ horo-ambient quasigeodesic in
$\widetilde{S}$ constructed from $\lambda^h \subset
{\widetilde{S}}^h$ \\ 
$\bullet$  $\beta^h$ $ = $ geodesic in $\widetilde{N^h}$ joining $a,
b$ \\
$\bullet$   $\beta_{amb}^0$  $ = $ horo-ambient quasigeodesic 
 in $\widetilde{N}$ obtained from  $\beta^h$ by replacement of
 hyperbolic by `Euclidean' geodesic segments for horoballs in
$\widetilde{N^h}$ \\
$\bullet$ $\beta_{amb} = \Pi_\lambda ( \beta_{amb}^0  ) \cap B_\lambda$ \\

\subsection{Proof of Theorem for Surfaces with Punctures}

By construction, the hyperbolic geodesic $\beta^h$ and the ambient
quasigeodesic  $\beta_{amb}^0$ agree exactly off
horoballs.  $\beta_{amb}$ is constructed from $\beta_{amb}^0$ by
projecting it onto $B_\lambda$ and so by Theorem \ref{mainref},
it is an ambient quasigeodesic. But it might
`backtrack'. Hence, we need to modify it such that it satisfies the no
backtracking condition. First, observe by Theorem \ref{ctm} that all
three $\beta^h$, $\beta_{amb}^0$, $\beta_{amb}$ track each other off
 some $K$-neighborhood of horoballs.

The advantage of working with $\beta_{amb}$ is that it lies on
$B_\lambda$.
However, it might backtrack.

\begin{lemma}
There exists $C > 0$ such that for all
 $x \in \lambda^b_i \subset B^b_\lambda \subset B_\lambda $
 if $\lambda^h$ lies outside
 $B_n(p)$ for a fixed reference point $p \in \widetilde{S^h}$ and also
 assuming that in fact the reference point lies in $\widetilde{S}$, then
 $x$ lies outside an $\frac{n-C}{C+1}$ ball about $p$ in
 $\widetilde{N}$.
\label{far1}
\end{lemma}

\noindent {\bf Proof:} 
Since $\lambda^{b}$ is a part of $\lambda^h$, 
therefore $r_x (0)$ lies outside $B_n (p)$. By Corollary
\ref{qgray2}, 
there exists $C > 0$ such that for all $i, j \in \{ 
0, 1, 2, \cdots$, 

\begin{center}

$|i-j| \leq d(r_x (i) , r_x (j)) \leq C|i-j|$

\end{center}

Also, $d(x, p) \geq i$ since $x \in \mu^b_i$. (Here distances are all
measured in $\widetilde{N}$ regarded as a tree of spaces with ``successive" vertex spaces separated by 
distance one.) Hence,

\begin{eqnarray*}
d(x,p) & \geq & min \{ i, n - C(i+1) \} \\
       & \geq & \frac{n-C}{C + 1}
\end{eqnarray*}

This proves the result. $\Box$

\smallskip

If $x \in B_\lambda  $, then $x \in B_\lambda^b$ or  $x \in
B_\lambda^c$.for some $\mu $.
 Hence $x \in B_\lambda $ implies that either $x$
lies on some horosphere bounding some $ \bf{H} \in \mathcal{H}$ or,
from Lemma \ref{far1} above, $d(x, p) \geq  \frac{n-C}{C + 1}$.
Since $\beta_{amb}$ lies on $B_\lambda $, we conclude that
$\beta_{amb}$ is an ambient quasigeodesic  in
$\widetilde{N}$ such that every point $x$ on $\beta_{amb}$ either
lies on some horosphere bounding some $ \bf{H} \in \mathcal{H}$ or,
from Lemma \ref{far1} above, $d(x, p) \geq  \frac{n-C}{C + 1}$.

McMullen \cite{ctm-locconn} shows (cf Theorem \ref{ctm} )
that in $\widetilde{N^h}$,
any such ambient quasigeodesic
$\beta_{amb}$ lies in a bounded neighborhood of $\beta^h \cup
{ \mathcal{H} } ( \beta^h )$.
We do not as yet know that $\beta_{amb}$ does not backtrack, but we
can convert it into one without much effort. 
Let $\Pi$ denote nearest point projection of  $\widetilde{N^h}$ onto
 $\beta^h \cup \mathcal{H}$ $ (\beta^h )$. 
Then $\Pi ( \beta_{amb} ) = \beta_1$ is again an
ambient quasigeodesic in $\widetilde{N}$. Further, $\beta_1$ tracks
$\beta_{amb}$ throughout its length, since $\Pi$ moves points through
a uniformly bounded distance. Now $\beta_1$ might backtrack, but it can
do so in a trivial way, i.e. if $\beta_1$ re-enters a horoball after
leaving it, it must do so at exactly the point where it leaves
it. Removing these `trivial backtracks', we obtain an  {\bf
  ambient quasigeodesic without backtracking} $\beta$ which tracks
$\beta_{amb}$ throughout its length. 

\smallskip

{\bf Note:} On the one hand $\beta$ is an ambient quasigeodesic
without backtracking. Hence, it reflects the intersection pattern of
$\beta^h$ with horoballs. On the other hand, it {\bf tracks}
$\beta_{amb}$
whose properties we already know from Corollary \ref{qgray2} above.

\smallskip

Since, of $\beta$ and $\beta^h$, one is an ambient quasigeodesic
without backtracking, and the other a hyperbolic geodesic joining the
same pair of points, we conclude from Theorem \ref{farb} that
they have similar intersection patterns with horoballs, i.e.
there exists $C_0$ such that \\
$\bullet$ If only of $\beta$ and $\beta^h$ penetrates or travels along the
  boundary of a horoball 
  $\bf{H}$, then it can do so for a distance $ \leq C_0$. \\
$\bullet$ If both $\beta$ and $\beta^h$ enter (or leave) a horoball
  $\bf{H}$ then their entry (or exit) points are at a distance of at
  most $C_0$ from each other.

Again, since $\beta$ tracks $\beta_{amb}$, we conclude that there
 exists
$C > 0$ such that $\beta$ lies in a $C$-neighborhood of $\beta_{amb}$
and hence from Lemma \ref{far1} \\
$\bullet$ Every point $x$ on $\beta$ either
lies on some horosphere bounding some $ \bf{H} \in \mathcal{H}$ or,
$d(x, p) \geq  \frac{n-C}{C + 1} -C$

\smallskip

The above three conditions on $\beta$ and $\beta^h$ allow us to 
deduce the following  condition for $\beta^h$.

\begin{prop} 
Let $x, \lambda^h$ satisfy the hypotheses of Lemma
\ref{far1}. 
Every point $x$ on $\beta^h$ either
lies inside some horoball $ \bf{H} \in \mathcal{H}$ or,
$d(x, p) \geq  \frac{n-C}{C + 1} -C = m(n)$
\label{far2}
\end{prop}

We split $\beta^h$ into two parts. $\beta^c$ consists of those points
of $\beta^h$ which lie within horoballs. We set $\beta^b$ to be the
closure of $\beta^h - \beta^c$. 

We have denoted $ \frac{n-C}{C + 1} -C$ by $ m(n)$, so that $m(n)
\rightarrow \infty$ as $n \rightarrow \infty$.
The above Proposition asserts that the geodesic $\beta^h$ lies outside
large balls about $p$ modulo horoballs. By Lemma \ref{contlemma} this
is almost enough to guarantee the existence of a Cannon-Thurston map. 
The rest of the necessary work is given below.

\begin{theorem}
Suppose $S^h$ is a hyperbolic surface of finite volume.
Suppose that $N^h$ is a hyperbolic manifold corresponding to a
 representation
of $\pi_1 (S^h )$ without accidental parabolics.
 Let $i : S^h \rightarrow N^h$ be a proper homotopy equivalence.
Then   $\tilde{i}: {\widetilde{S}}^h \rightarrow \widetilde{N^h}$ 
extends continuously to the boundary
$\hat{i}: \widehat{S^h}
\rightarrow \widehat{N^h}$. If $\Lambda$ denotes the limit set of
$\widetilde{M}$, then $\Lambda$ is locally connected.
\label{main-punct}
\end{theorem}

\noindent {\bf Proof:} 
Let $\lambda^h$ be a geodesic segment in 
$\widetilde{S^h}$
lying outside $B_n{(p)}$ for some fixed reference point $p$. Fix 
neighborhoods of the cusps and lift them to the universal cover. Let
$\mathcal{H}$ denote the set of horoballs. Assume without loss of
generality that $p$ lies outside horoballs. Let $\beta^h$ be 
the hyperbolic geodesic in $\widetilde{N^h}$ joining the endpoints of
$\lambda^h$. Further, let $\beta^h = \beta^b \cup \beta^c$ as above.
Then by Proposition \ref{far2}, $\beta^b$ lies outside an $m(n)$ ball
about $p$, with 
 $m(n)
\rightarrow \infty$ as $n \rightarrow \infty$.

Next, let ${\bf{H}}_1$ be any one of the equivariant collection of horoballs 
that $\beta^h$ meets. By shrinking the cusps
in the quotient manifold slightly if necessary, we may assume without loss of generality
that none of the equivariant collection of horoballs  contain $p$. The
 entry and exit
points $u$ and $v$ of $\beta^h$ into and out of  ${\bf{H}}_1$
lie  outside an $m(n)$ ball
about $p$. Let $z$ be the point on the boundary sphere that
${\bf{H}}_1$ is based at. Then for any sequence  $x_i \in
{\bf{H}}_1$
with $d(p, x_i ) \rightarrow \infty$, it follows that
$x_i \rightarrow z$. 
Let  $\{ x_i \}$ and $\{ y_i \}$ denote two such sequences on any horoball ${\bf{H}}_1$ not 
containing $p$. 
Then the
visual diameter of the set $\{ x_i , y_i \}$ must go to zero. Hence,
if $[x_i,y_i]$ denotes the geodesic joining $x_i , y_i$ then
$d(p,[x_i,y_i]) \rightarrow \infty$. Since, $u, v$ lie outside an
$m(n)$ ball, there exists some function $\psi$
(independent of ${\bf{H}}_1$), such that the geodesic
$[u,v]$
lies outside a $\psi (m(n))$ ball around $p$, where $\psi (k)
\rightarrow \infty$ as $k \rightarrow \infty$. 

Since the choice of this function does not depend on
 ${\bf{H}}_1$, which  is chosen at random, we conclude that there
exists such a function for all of $\beta^c$. We have thus established:
 \\
$\bullet$ $\beta^b$ lies outside an $m(n)$ ball about $p$. \\
$\bullet$ $\beta^c$ lies outside a $\psi (m(n))$ ball about $p$. \\
$\bullet$ $ m(n)$ and  $\psi (m(n))$ tend to infinity as $n \rightarrow
 \infty$\\

Define $f(n) = min( m(n),\psi (m(n)) )$. 
Then  $\beta^h$ lies outside an $f(n)$ ball about $p$ and 
 $f(n)
\rightarrow \infty$ as $n \rightarrow \infty$.

 By Lemma \ref{contlemma}  $i: \widetilde{S^h}
\rightarrow \widetilde{N^h}$ extends continuously to the boundary
$\hat{i}: \widehat{S^h}
\rightarrow \widehat{N^h}$. This proves the first statement of the
theorem.

Now, the limit set of ${\widetilde{S}}^h$  is the circle at infinity,
which is
locally connected.  Further, the
continuous image of a compact locally connected set is locally
connected \cite{hock-young}. Hence, if $\Lambda$ denotes the limit set of
$\widetilde{N^h}$, then $\Lambda$ is locally connected. This proves
the theorem. $\Box$

\bibliography{bddgeo}
\bibliographystyle{alpha}

\noindent RKM Vivekananda University, 
Belur Math, India \\
E-mail: mahan@@rkmvu.ac.in

\end{document}